\documentstyle[12pt]{report}

\newcommand{\vs}{\vspace{.15in}}

\newcommand{\noin}{\noindent}

\newcommand{\prop}{{\bf Proposition.} }
\newcommand{\lemma}{{\bf Lemma.}  }
\newcommand{\proof}{{\sl Proof.}  }
\newcommand{\rem}{{\bf Remark.} }
\newcommand{\defin}{{\bf Definition.} }
\newcommand{\cor}{{\bf Corollary.} }

\newcommand{\qed}{{$\Box $} \vs}
\def\C{{\bf C}}

\def\hZ{{\hat Z}}

\def\hJ{{\hat J}}

\def\tY{{\tilde Y}}

\begin{document}
\begin{center} {\bf SUBRINGS INVARIANT UNDER ENDOMORPHISMS}
\\[3ex] {\bf SHULIM KALIMAN} \footnote{ The author was partially
supported by the NSA grant MDA904-00-1-0016.}\\[3ex] Department of
Mathematics, University of Miami \\ Coral Gables, Florida 33124,
USA \\[4ex]
\end{center}
Abstract. {\small Let $S$ and $R$ be the rings of regular
functions on affine algebraic varieties over a field of
characteristic 0, $R$ be embedded as a subring in $S$, and $F : S
\to S$ be an endomorphism such that $F(R) \subset R$. Suppose that
every ideal of height 1 in $R$ generates a proper ideal in $S$,
and the spectrum of $R$ has no selfintersection points. We show
that if $F$ is an automorphism so is $F|_R : R \to R$. When $R$
and $S$ have the same transcendence degree then the fact that
$F|_R$ is an automorphisms implies that $F$ is an automorphism.}
\\[3ex]

{\bf 1. Introduction.} In [CZ] E. Connell and J. Zweibel proved
the following fact. Let $k$ be a field of characteristic 0, $S$
and $R$ be isomorphic to $k[x_1, \ldots , x_n]$, $R$ be a subring
of $S$, and $F: S \to S$ be an endomorphism for which $F(R)
\subset R$. Then $F$ is an automorphism iff $F|_R : R \to R$ is an
automorphism.

Though the result is very natural the proof is not simple and it
is based to a great extend on the Zariski Main theorem. We shall
study the question when an analogue of this theorem holds for a
wider class of rings. If one suppose that $S$ and $R \subset S$
are the rings of regular functions on affine algebraic varieties
(over $k$) then a similar theorem is not valid without an extra
assumption. Put $S = k[x,x^{-1},y]$ and $R=k[x,y]$. Consider the
automorphism of $S$ that sends $x,x^{-1}, y$ to $x, x^{-1}, xy$
respectively. Then its restriction to $R$ is not an automorphism
though the image of $R$ is contained in $R$. This counterexample
is based on the fact that $x$ is a unit in $S$ but not in $R$.
Meanwhile it is easy to check that under the assumption of the
Connell-Zweibel theorem every element of the subring which is
invertible in the ambient ring must be automatically invertible in
the subring. It turns out that this property is crucial in the
case when $R$ is a UFD. In a more general setting we prove \vs

{\bf Theorem A.} {\sl Let $S$ and $R$ be affine domains over a
field $k$ of characteristic 0, $R$ be embedded as a subring in
$S$, and $F : S \to S$ be an endomorphism for which $F(R) \subset
R$.

(i) Suppose that $R$ is the ring of regular functions on an affine
algebraic variety without selfintersection points \footnote{ The
absence of selfintersection points is essential. Indeed, in the
example above we can replace $R$ by its subring which consists of
polynomials in $k [x,y]$ taking the same values at points (0,0)
and (1,1). This gives a counterexample to Theorem A in the
presence of selfintersection points.} (for instance, $R$ is
integrally closed) and every ideal of height 1 in $R$ generates a
proper ideal in $S$. Then if $F$ is an automorphism so is $F|_R :
R \to R$.

(ii) Let $S$ and $R$ have the same transcendence degree.
Then if $F|_R$ is an automorphism so is $F$.} \vs

Using the ``Lefschetz principle" (e.g., see [BCW]) one can reduce
the problem to the case when $k = \C$. Furthermore, we prefer to
work with a geometrical reformulation of this theorem. More
precisely, Theorem A is a consequence of \vs

{\bf Theorem B.} {\sl Let $X$ and $Y$ be irreducible complex
affine algebraic varieties. Suppose that $\rho : X \to Y$, $f: X
\to X$, and $g: Y \to Y$ are morphisms such that $\rho$ is
dominant and the following diagram is commutative \vs

\[ \begin{array}{cccc}
X & \stackrel{f}{\rightarrow} & X & \,\\
\, \, \downarrow \rho && \, \, \downarrow \rho & \, \, \, \, \, \,
\, \, \, \, \, \, \, \, \, \, \, \, \, \, \, \, \, \,
\, \, \, \, \, \, \, \, \, \, \, \, (1) \\
Y & \stackrel{g}{\rightarrow} & Y & \,\\
\end{array} \] \vs

(i) Suppose that $Y$ has no selfintersection points, $f$ is an
automorphism, and $g$ is not. Then there exists a closed
hypersurface $D \subset Y$ such that ${\rm codim }_Y g(D) \geq 2$
and $\rho^{-1} (D)$ is empty.

(ii) Let $dim \, Y = \dim X$. Then if $g$ is an automorphism so is
$f$.} \vs

Besides the Zariski Main Theorem [H, Ch. 5, Th. 5.2] our other
main tool follows from a remarkable theorem of Ax [A] (later
rediscovered by Kawamata [I]) \vs

{\bf Theorem.} {\sl Let  $Z$ be a complex algebraic variety and
let $h : Z \to Z$ be an injective morphism. Then $h$ is an
automorphism.} \vs

The idea of the proof is the following. Using the Zariski and Ax
theorems we prove that if $g$ (resp. $f$) is not an automorphism
under the assumption of Theorem B (i) (resp. B(ii)) then there
exists a divisor $D \subset Y$ (resp. $E \subset X$) such that
${\rm codim}_Y g(D) \geq 2$ and $g(D)\subset D$ (resp. ${\rm
codim}_X f(E) \geq 2$ and $f(E) \subset E$). The next argument is
especially simple in the smooth equidimensional  case: we show
that the zero multiplicity of the Jacobians of $g^s \circ \rho$
and $\rho \circ f^s$ are different at $x \in \rho^{-1}(D)$ (resp.
$x \in E$) for some $s>0$. In the non-smooth case we show that the
dimensions of the images of a $k$-jet space at $x$ under $g^s
\circ \rho$ and $\rho \circ f^s$ are different.

It is our pleasure to thank M. Miyanishi for drawing our attention
to the paper of Ax. \vs

{\bf 2. The existence of the exceptional divisor.\\ 2.1.}
Replacing $X$ and $Y$ in diagram (1) with their normalizations
$X^0$ and $Y^0$ (which are also affine) we get a commutative
diagram \vs

\[ \begin{array}{ccc}
X^0 & \stackrel{f^0}{\rightarrow} & X^0\\ \, \, \downarrow \rho^0
&& \, \, \downarrow \rho^0 \\ Y^0 & \stackrel{g^0}{\rightarrow} &
Y^0\\
\end{array} \] \vs

As $Y$ has no selfintersection points the normalization $Y^0 \to
Y$ is a homeomorphism. Hence for any divisor $D \subset Y$ and its
proper transform $D^0\subset Y^0$ we have $(\rho^0) ^{-1}(D^0) \ne
\emptyset$ iff $\rho^{-1}(D) \ne \emptyset$. Hence it is not
difficult to prove the following.

 \lemma {\sl Theorem B is true if
it is true under the additional assumption that $X$ and $Y$ are
normal.} \vs

%
%

{\bf 2.2.}
\lemma {\sl Let $X$ and $Y$ be as in diagram (1). Then

(a) if $f$ is birational so is $g$,

(b) if $dim \, X = dim \, Y$ and $g$ is birational
then $f$ is birational.}

{\sl Proof.} Consider (a). It follows from the semi-continuity
theorem [H, Ch. 3, Th. 12.8] that the number of connected
components in $\rho^{-1} (y)$ is an upper semi-continuous function
on $Y$. In particular, this number is the same for general points
$y \in Y$. Denote it by $n$. Note that $g$ is dominant since
otherwise $f$ is not dominant. Let $k$ be the number of components
in the preimage of a general point of $y \in Y$ under $g$. There
are $n$ components in $(\rho \circ f)^{-1}(y)$ and $k n$
components in $(g \circ \rho) ^{-1}(y)$ . By commutativity of
diagram (1) we have $k=1$. That is, the degree of $g$ is 1 and $g$
is birational. The proof of (b) is similar. \qed

\cor {\sl Under the assumption of Theorem B $f$ is birational iff
$g$ is birational.}\\

{\bf 2.3.} By the semi-continuity theorem  $X_1 = \{ x \in X |dim
f^{-1} (f(x)) > \dim X - \dim Y \}$ is a closed algebraic
subvariety of $X$. Let $X_0 = X \setminus X_{1}$ and $Y^0$ be the
largest Zariski open subset of $\rho (X_0)$. In Theorem B (ii) we
need also the Zariski open subset $X^0$ of $X$ that is the largest
subset such that $\rho |_{X^0}$ is quasi-finite.

\lemma {\sl (1) Under the assumption of Theorem B (i) the
restriction of $g$ to $Y^0$ is an automorphism provided that $Y$
is normal.

(2) Under the assumption of Theorem B (ii) the restriction of $f$
to $X^0$ is an automorphism provided that $X$ is normal. }

{\sl Proof.} The the commutativity of diagram (1) implies that
$f(X_1) \subset X_1$ in the first statement. By the Ax theorem the
restriction of $f$ to $X_1$ is an automorphism of $X_1$ whence we
have the similar fact for $X_0$. The commutativity of diagram (1)
implies that the restriction of $g$ to $\rho (X_0)$ is a
homeomorphism of $\rho (X_0)$ whence (1) follows from the Zariski
Main theorem. In (2) let $E \subset X$ be the set of points where
$f$ is not \'etale. By the Zariski Main Theorem any $x \in E$ is
not a connected component of $f^{-1}(f(x))$, and by the
commutativity of diagram (1) $f^{-1}(f(x))$ is contained in
$\rho^{-1}(\rho (x))$. Thus $X^0 \subset X\setminus E$. As $\rho
=g\circ \rho \circ f^{-1}$ for $x \in X^0$ we have $f(x)\in X^0$,
i.e. $X^0 \subset X^0$ whence by the Ax theorem $f|_{X^0}: X^0 \to
X^0$ is an automorphism. \qed

{\bf 2.4.}
\prop {\sl Let $g: Y \to Y$ be a birational endomorphism of a
normal affine algebraic variety which is not an autmorphism, but
for a Zariski open subset $Y^0$ of $Y$ the restriction of $g$ to
$Y^0$ is an automorphism. Then there exists an exceptional divisor
$D$ with respect to $g$ (i.e. ${\rm codim}_Y g(D) \geq 2$).
Furthermore, replacing $g$ with $g^m$ for some $m>0$ one can
suppose that $g(D) \subset D$. }

{\sl Proof.} Let $D' =Y \setminus Y^0$.  Denote by $D_0'$ the
Zariski open subset of $D'$ that consists of points such that the
restriction of $g$ to a neighborhood of any of these points is a
quasi-finite morphism. For every $y\in D'$ its image $g(y)$ cannot
belong to $Y^0$ (that is, $g(D') \subset D'$) since otherwise
$g(y_1) = g(y)$ for some $y_1 \in Y^0$ whence the preimage of
$g(y)$ is not connected contrary to the Zariski Main theorem. The
same theorem implies that the restriction of $g$ to $Y^0 \bigcup
D_0'$ is an embedding. Suppose that $C$ is an irreducible
component of $D'$ which is a hypersurface and which meets $D_0'$
(i.e. $D_0'\cap C$ is dense in $C$), and let $D_1$ be the union of
such hypersurfaces. Then the closure of $g(C)$ is also a
hypersurface which is an irreducible component of $D'$. Assume
that this component is not contained in $D_1$. Denote by $D_0''$
the subset of $D'$ that consists of points such that the
restriction of $g^2$ to a neighborhood of any of these points is a
quasi-finite morphism. Note that under this assumption $C$ does
not meet $D_0''$. Thus replacing, if necessary, $g$ with $g^m$ for
some natural $m$ we can suppose that $g(D_1) \subset D_1$. In
particular, $D_1\setminus g(D_1 \cap D_0')$ is of codimension at
least 2 in $Y$. Assume that $D'$ does not contain an exceptional
divisor with respect to $g$. Then the codimension of the
complement to $g(Y^0 \cup D_0')$ in $Y$ is at least 2. Since
$g^{-1}$ is well-defined on $g(Y^0 \cup D_0')$ it can be extended
to $Y$ by the theorem about deleting singularities for normal
algebraic varieties in codimension 2 [D, Ch. 7.1]. This
contradicts the assumption that $g$ is not an automorphism whence
there exists an exceptional divisor $D$ with respect to $g$ which
is, of course, contained in $D'$.

For the second statement note that for every $y \in D$ its image
$y_1 =g(y)$ must belong to an irreducible component of $D'$ which
is a hypersurface since otherwise $g^{-1}$ can be extended to
$y_1$ by the theorem about deleting singularities in codimension
2. Suppose that $C$ and $D_1$ are as above. In particular, the
closure of $g(D_1)$ is $D_1$, and $g(C \cap D_0')$ is dense in
$g(C)$. Let $C_0$ be the complement in the closure of $g(C)$ to
the union of the other components of $D'$ that are hypersurfaces.
Note that $g^{-1}(C_0)$ is contained in $D_0'$ by the theorem
about deleting singularities. Furthermore, applying this theorem
again we see that that $g^{-1}(D_1 \setminus D)$ is also contained
in $D_0'$, i.e. $y_1$ cannot belong to $D_1 \setminus D$. Thus
$y_1 \in D$ and $g(D) \subset D$. \qed

\cor {\sl If $g$ (resp. $f$) is not an automorphism under the
assumption of Theorem B (i) (resp. B(ii)) then there exists an
exceptional divisor $D$ with respect to $g$ (resp. $E$ with
respect to $f$). Furthermore, one can suppose that $g(D) \subset
D$ (resp. $f(E) \subset E$). } \vs

{\bf 2.5.} We can already prove Theorem B in the case of smooth
varieties $X$ and $Y$ (for simplicity we shall consider the case
when $X$ and $Y$ are of the same dimension). Consider a
holomorphic mapping $h : V \to U$ of equidimensional complex
manifolds $V$ and $U$ and the Jacobian of this mapping in local
coordinate systems at $v \in V$ and $u=h(v)$, i.e. the determinant
of the Jacobi matrix. The Jacobian itself depends on the choice of
these local coordinate systems but the order of its zeros at $v$
does not. We denote this order by $Jd_h(v)$ . The following the
two properties of $Jd_h$ are simple.

($\alpha$) $Jd_h(v) >0$ iff $h$ is not a local embedding in a
neighborhood of $v$;

($\beta$) if $e: U \to W$ is another holomorphic mapping of
equidimensional complex manifolds then $Jd_{e \circ h} \geq Jd_h
(v) + Jd_e (u)$, and the equality holds in the case when either
$h$ is a local embedding at $v$ or $e$ is a local embedding at
$u$.

Let the assumption of Theorem B (i) hold and $D$ be as in
Corollary 2.4. Assume that $\rho^{-1} (D) \ne \emptyset$, $x \in
\rho^{-1} (D), x'=f(x)$, and $y =\rho (x)$. Since $\rho \circ f =
g \circ \rho$ and $f$ is an automorphism  we have by ($\beta$)
$Jd_{\rho \circ f}(x) = Jd_{\rho} (x') = Jd_{\rho}(x) + Jd_g (y)$.
Since $g$ is not a local embedding at $y$ we see that $Jd_g(y)
>0$. Furthermore, since $g(D) \subset D$, replacing
$g$ (resp. $f$) by $g^m$ (resp. $f^m$) we can make $Jd_g(y) >>0$.
One the other hand $Jd_{\rho} (x')$ is bounded as $Jd_{\rho}$ is
bounded on $X$. This contradiction concludes the proof of Theorem
B (i) in the smooth case.
The proof of Theorem B (ii) in the smooth case is similar. \vs

{\bf 3. Jets on manifolds.\\ 3.1.} In order to deal with the
general case we need to consider the variety of $k$-jets $J^k(M)$
from the germ $(\C , 0)$ of the complex line at the origin into a
complex manifold $M$. The following notation and simple facts will
be used. For $k \geq l$ we denote by $\tau_{M}^{k,l} : J^k(M) \to
J^l(M)$ the natural projection. The map $\tau_M^{k,0} : J^k(M) \to
J^0(M) \simeq M$ is a $\C^s$-fibration where $s= k \dim M$. This
fibration admits a natural $\C^*$-action generated by the
$\C^*$-action on $(\C , 0)$. The restriction of this action to any
fiber generates an embedding of this fiber into a weighted
projective space. Hence we can extend $\tau_M^{k,0}$ to a proper
holomorphic fibration $ \bar{\tau}^{k,0}_M: \bar{J}^k(M) \to M$
whose fibers are isomorphic to this weighted projective space. For
every subset $Z$ of $J^l(M)$ we denote by $J_Z^k(M)$ the set $\{ j
\in J^k(M) \, | \, \tau_M^{k,l}(j) \in Z \}$. Note that if $Z$ is
a variety then $$\dim J_Z^k(M) = \dim Z + (k-l) \dim M. \eqno
(2)$$ Any holomorphic map of complex manifolds $\varphi : M \to N$
generates a holomorphic map $\varphi^{(k)} : J^k(M) \to J^k(N)$
such that $\tau_N^{k,l} \circ \varphi^{(k)} = \varphi^{(l)} \circ
\tau_M^{k,l}$. In particular, if $Z \subset J^l(M), Z^k \subset
J_Z^k(M)$, $W=\varphi^{(l)} (Z)$, and $W^k =\varphi^{(k)} (Z^k)$
then $W^k \subset J_W^k (N)$. Another useful observation is that
$\varphi^{(k)} $ commutes naturally with the $\C^*$-actions on
$J^k(M)$ and $J^k(N)$ whence it can be extended to a holomorphic
map $ \bar{\varphi}^{(k)} : \bar{J}^k (M) \to \bar{J}^k (N)$. \vs

{\bf 3.2.} \prop {\sl Let $\varphi : M \to N$ be a non-degenerate
holomorphic map of complex manifolds. Let $l \geq 0$ and $Z_0$ be
an algebraic subvariety of $J^l(M)$. Then there exists $r \geq l$
such that for every $k \geq r$, $Z=J_{Z_0}^r(M), Z^k =J_Z^k(M)$,
$W=\varphi^{(r)}(Z)$, and $W^k =\varphi^{(k)}(Z^k)$ we have $\dim
W^k =\dim W + (k-r) \dim N$ (i.e., by (2), $W^k$ is dense in
$J_{W}^k(N)$).}

\proof First note that we can suppose that $Z_0$ is irreducible
(then in the non-irreducible case for every irreducible component
of $Z_0$ one can find its own number $r$ and take the maximum of
these numbers in the statement of Lemma).

Step 1. Let us show that for every $l \geq 0$ it suffices to prove
the statement under the additional assumption that $\tau^{l,0}_M
(Z_0)$ is a point $x_0 \in M$.

Let $B=\tau^{r,0}_M(Z)= \tau_M^{l,0} (Z_0)$ and $x \in B$. Put
$\theta^k = \tau_N^{k,r} |_{W^k} : W^k \to W$, $Z_x^k = \{ j \in
Z^k | \, \tau^{r,0}_M(j) =x \}$, and $W_x^k =\varphi^{(r)}
(Z_x^k)$. Since $\theta^k (W_x^k) =W_x^r$ it suffices to show that
for general $x \in B$ we have $\dim W_x^k -\dim W_x^r =(k-r) \dim
N$. Since we assume that Lemma is correct under the additional
assumption, for every $x \in B$ the number $r$ can chosen so that
we have this equality. Furthermore, by Baire's category theorem we
can suppose that there exists $r$ for which the equality holds for
every $x$ in a subset $L \subset B$ which is not contained in any
analytic subset of $B$. Consider $\tilde{W}^k$ equal to the image
of $Z^k$ in $B \times W^k$ under the holomorphic map $(
\tau^{k,0}_M, \varphi^{(k)} )$. Note that $W_x^k$ can be viewed as
a fiber of the natural projection $ \tilde{W}^k \to \tau_M^{k,0}
(Z^k)=B$. As $\varphi^{(k)}$ and $\tau_M^{k,0}$ can be extended to
$ \bar{\varphi}^{(k)}$ and $ \bar{\tau}_M^{k,0}$ from 3.1 this
projection can be extended to a proper holomorphic map into $B$
whence, by semi-continuity theorem (e.g., see [BN, Th. 2.3]) the
dimension of $W_x^k$ is constant on a complement $U$ to a proper
analytic subset of $B$. Consider the natural projection
$\tilde{\theta}^k : \tilde{W}^k \to \tilde{W}^r$ generated by
$\tau_N^{k,r}$ and its restriction to $W_x^k$ which may be viewed
as $\theta^k$. The dimension of a general fiber of
$\tilde{\theta}^k$ is $\dim \tilde{W}^k - \dim \tilde{W}^r$ where
the last number coincides with $\dim W_x^k -\dim W_x^r$ for
general $x$, i.e. for $x \in U$. Since $L$ meets $U$ we see that
$\dim W_x^k -\dim W_x^r =(k-r) \dim N$ for $x \in U$ which
concludes the first step.

{\sl Step 2.} If $Z \subset J_{x_0}(M)$ the fact becomes local
analytic, and one can suppose that $M$ (resp. $N$) coincides with
the germ $(\C^m , o_m)$ (resp. $(\C^n , o_n) $) of a Euclidean
space at the origin (of course, we put $x_0 =o_m$). Let
$(\varphi_1, \ldots , \varphi_n)$ be the coordinate form of a
holomorphic map $\varphi$ and let $\varphi_{i, 0}$ the the minor
homogeneous form in the Taylor decomposition of $\varphi_i$. We
need

{\bf Claim.} {\sl For $Z_0 \subset J_{x_0} (M)$ it suffices to
prove the local version of Lemma in the case of homogeoneous
$\varphi$, i.e. $\varphi_i =\varphi_{i,0}$ for every $i$ and the
degrees of these coordinate functions are the same number $s$.}

First note that if $\theta : \C^n \to \C^n$ is a polynomial map
Lemma holds for morphism $\varphi$ provided it holds for morphism
$\phi = \theta \circ \varphi$. The coordinate functions $\phi_1,
\ldots , \phi_n$ of $\phi$ are elements of the algebra generated
by $\varphi_1, \ldots , \varphi_n$. These elements can be chosen
so that there minor homogeneous forms are algebraically
independent [M-L]. Thus we can suppose from the beginning that
$\varphi_{1,0}, \ldots , \varphi_{n,0}$ are algebraically
independent (i.e. morphism $\psi_0 = (\varphi_{1,0}, \ldots ,
\varphi_{n,0})$ is dominant). Furthermore, replacing $\varphi_1 ,
\ldots , \varphi_n$ by their powers, we suppose that each
$\varphi_{i,0}$ has the same degree $s$. Let $ { \xi} = (\xi_1 ,
\ldots , \xi_m)$ be a coordinate system on $\C^m$. Put $\psi_{i,c}
=c^{-s} \varphi_{i,0} (c { \xi})$ where $c \in \C^*$, and put
$\psi_c = ( \psi_{1,c}, \ldots , \psi_{n,c})$. Clearly,
$\varphi^{(k)} (J_{Z}^k (M))$ and $\psi_{c}^{(k)} (J_{Z}^k(M))$
are isomorphic for $c \ne 0$, and $\psi_c \to \psi_0$ as $c \to
0$. This yields a surjective morphism from $\varphi^{(k)} (J_{Z}^k
(M))$ to $\psi_{0}^{(k)} (J_{Z}^k(M))$ which implies the statement
of the Claim and concludes Step 2.

Step 3. We shall use induction by $l$. Let $l=0$. By Step 1 we can
suppose that $Z_0$ consists of one element $j_0$ which is
presented by a constant map from $(\C , 0)$ into a point $x_0 \in
M$. That is, $j_0(t) =x_0$ where $t$ is a coordinate on $(\C ,0)$.
By Step 2 we can suppose that $M=\C^m, x_0 =o_m, N=\C^n$, and
$\varphi : \C^m \to \C^n$ is homogeneous of degree $s$. If $n=m$
then, since $\varphi$ is dominant, it is a local analytic
isomorphism at a general point $x$ of $\C^n$. Hence for $y=
\varphi (x)$ the restriction of $\varphi^k$ to $J_x^k(M)$ is an
isomorphism between $J_x^k(M)$ and $J_y^k(N)$. In the case when
$m>n$ applying the above argument to the restrictions of $\varphi$
to general $n$-dimensional submanifolds of $M$ we can see that the
restriction of $\varphi^{(k)}$ to $J_{x}^k(M)$ is an epimorphism
onto $J_{y}^k(N)$ for general $x \in M$. Every $j \in
J_{o_m}^k(M)$ is of form $$j(t) =tj_1(t) \eqno (3)$$ where $j_1
\in J_{x}^{k-1}(M)$ and $x \in \C^n$. Put $r=s$ and consider the
Zariski open subset of $J_{o_m}^r(M)$ which consists of $j^0$ such
that $j^0(t) =tj_1^0(t)$ where $j_1^0$ is an element of
$J^{r-1}(M)$ for which $x=j_1^0(0)$ is a general point of $\C^n$.
In particular, $j^0$ is a general element of $Z=J_{o_m}^r(M)$, and
the restriction of $\varphi^{(k)}$ to $J_x^k(M)$ is an epimorphism
onto $J_y^k(N)$ where $y =\varphi (x)$. Let $j \in J_{j^0}^k(M)$
and $j_1 \in J_{j_1^0}^{k-1}$ be as in (3). Note that
$$\varphi^{(k)} (j) = t^s \varphi^{(k-s)}(j_2) \eqno (4)$$ where
$j_2= \tau_M^{k-1,k-s} ( j_1)$. Hence $\varphi^{(r)} (j^0) = t^s
y$. Since the restriction of $\varphi^{(k-s)}$ to $J_x^{k-s}(M)$
is an epimorphism onto $J_y^{k-s}(N)$ we see that the restriction
of $\varphi^{(k)}$ to $J_{j^0}^k(M)$ is an epimorphism onto
$J_{\varphi (j^0)}^k(N)$ which proves the statement for $l=0$ and
concludes Step 3.

Step 4. Assume that Lemma is proven for $l-1$. That is, for every
$Z_0' \subset J_x^{l-1} (M)$ there exists $r_0 \geq l-1$ such that
for $Z'= J_{Z_0'}^{r_0}(M)$, $W'= \varphi^{(r_0)}(Z')$, and every
$k \geq r_0$ the image $\varphi^{(k)} (J_{Z'}^k (M))$ is dense in
$J_{W'}^k (N)$. By Step 1 we can suppose that
$\tau^{l,0}_M(Z_0)=x_0$ whence by Step 2 $M=\C^m, x_0 =o_m,
N=\C^n$, and $\varphi : \C^m \to \C^n$ is homogeneous of degree
$s$. This means that $Z_0$ is of form $Z_0 =t Z_0'$ and $Z=tZ'$.
Put $r=r_0+s$, $Z'' = \tau_M^{r-1,r-s} (Z')$ and $W'' =
\varphi^{(r-s)}(Z'')$. Then $Z'' =J_{Z_0'}^{r_0}(M)$ since
$\tau_M^{r-1,l-1}= \tau_M^{r-s,l-1} \circ \tau_M^{r-1,r-s}$. By
(4), $W=t^s W''$ and the statement of Lemma is equivalent to the
fact that $\varphi^{(k-s)}(J_{Z''}^{k-s}(M))$ is dense in
$J_{W''}^{k-s}(N)$. But this is true by the induction assumption
for $l-1$. \qed

{\bf 4. Jets on algebraic varieties. \\ 4.1.} We need an analogue
of $J^k(M)$ in the case of non-smooth algebraic varieties. In the
rest of the paper for every algebraic variety (resp. analytic set)
$Y$ and $y \in Y$ we denote by $(Y,y)$ the germ of $Y$ at $y$ in
the Zariski (resp. Euclidean) topology. Let $(Y,y) \hookrightarrow
(\C^n , o_n)$ be a closed embedding where $o_n$ is the origin in
$\C^n$. Let $t$ be a coordinate on $(\C ,0)$. We denote by $
\hat{J} \C^n$ the set of formal jets $ \hat{j}$ which are
$n$-tuples $ \hat{j} = ( \hat{j}_1, \ldots , \hat{j}_n)$ of formal
power series in $t$. Its subset $ \hat{J}_{o_n}\C^n$ consists of $
\hat{j}$ such that $ \hat{j}_i(0)=0$ for every $i$. We define the
set of formal jets $ \hat{J}_yY$ of $Y$ at $y$ as a subset of $
\hat{J}_{o_n} \C^n$ such that $ \hat{j} \in \hat{J}_yY$ iff for
every regular function $h$ from the defining ideal of $(Y,y)$ in
$(\C^n, o_n)$ the formal series $h \circ \hat{j}$ is zero. \vs

\defin Let $\tau^{k} : \hat{J} \C^n \to J^k\C^n$ be the forgetting projection.
The set of $k$-jets of $Y$ at $y$ is $J^k_yY : = \tau^k
(\hJ^k_yY)$. \vs

\rem We call $ \hat{j} \in (\tau^{k}) ^{-1}(j)$ a formal extension
of $j\in J^k\C^n$. By Artin's theorem [P, Th. 4.4] for $j \in
J^k_yY$ its formal extension $ \hat{j} \in \hat{J}_yY$ can be
chosen convergent. That is, we can treat $ \hat{j}(t)$ as a germ
of a curve in $Y$.\vs

{\bf 4.2.} \lemma {\sl The closure of $J_y^kY$ in $J_{o_n}^k \C^n$
is an algebraic variety, it is independent (up to an isomorphism)
from the choice of a coordinate $t$ on $(\C , 0)$ and from the
choice of the closed embedding $(Y,y) \hookrightarrow (\C^n ,
o_n)$, and $\tau_{Y}^{k,l}(J_y^kY) =J_y^lY$ where
$\tau_Y^{k,l}=\tau_{\C^n}^{k,l}|_{J_y^{k}Y}$ and $l \leq k$.
Furthermore, any morphism $\varphi : (Y,y) \to (Z,z)$ generates a
morphism $\varphi^{(k)} : J_y^kY \to J_z^kZ$. }

\proof For the first statement note that $ \hat{J}_yY$ is given in
$ \hat{J}_{o_n} \C^n$ by a countable number of polynomial
equations on the coefficients of the coordinates $ \hat{j}_i$ of
formal jets $ \hat{j}= ( \hat{j}_1, \ldots , \hat{j}_n)$. This
implies that $ \tau_{\C^n}^{k} ( \hat{J}_yY)$ is the intersection
of at most countable number of constructive sets whence the
closure of $ {J}_y^kY$ is an algebraic variety. The other
statements are immediate consequence of the definition. \qed

{\bf 4.3.} Consider a coordinate form $(j_1(t), \ldots , j_n(t))$
of a $j\in J_{o_n}^k\C^n$ where $t \in (\C , 0)$ and each $j_i$ is
a polynomial in $t$ of degree at most $k$. We say that the
multiplicity of $j$ is $m = \min \{ s | \exists l: {\frac
{d^s}{dt^s}} j_l (0) \ne 0 \}$. The subset of jets of multiplicity
$m$ in $J_{o_n}^k(\C^n)$ will be denoted by $J_{o_n}^{k,m}(\C^n)$.
Notation $ \hat{J}_{o_n}^m\C^n, \hat{J}_y^mY, J_y^{k,m}Y$ have the
similar meaning. For any $h$ from the defining ideal of $(Y,y)$ in
$(\C^n, o_n)$ consider its homogeneous decomposition $h=h_0+h_1+
\ldots$ where $h_0$ is the minor homogeneous form. One can treat
the tangent space of $\C^n$ at $o_n$ as 1-jets. Then the reduced
tangent cone $C_yY$ consists of all 1-jets $j(t)$ such that $h_0
\circ j (t)=0$ for any $h_0$ as above. This implies. \vs

\lemma {\em Every $j \in J^{k,k}_yY$ is of form $j(t) = j^1(t^k)$
where $j^1 \in C_yY$.} \vs

{\bf 4.4.} Let $ \sigma : \tilde{\C}^n \to (\C^n, o_n)$ be the
blowing-up of $(\C^n, o_n)$ at $o_n$ and $E$ be its exceptional
divisor. Let $(\xi_1, \ldots , \xi_n)$ be a coordinate system on
$(\C^n , o_n)$. Then $ \tilde{\xi} =( \tilde{\xi}_1 , \ldots ,
\tilde{\xi}_n)= (\xi_1, \xi_2/\xi_1, \ldots , \xi_n/\xi_1)$ is a
local coordinate system on $ \tilde{\C}^n$. Without loss of
generality we can suppose that ${ \frac {d^m}{dt^m}} j_1(0) \ne
0$. Then $(j_1, j_2/j_1, \ldots , j_n/j_1)$ can be viewed as an
$n$-tuple of power series so that the first $(k-m)$ terms of every
entry are well-defined. This enables us to define for $l =0,
\ldots , k-m$ morphism $ \theta_n^{k,m,l} : J_{o_n}^{k,m}\C^n \to
J_{E}^l \tilde{\C}^n$ such that in this local coordinate system
$\tilde{\xi}$ we have $\theta_n^{k,m,l}(j) = ([j_1]_{l},
[j_2/j_1]_{l}, \ldots , [j_n/j_1]_{l})$ where for every power
series $a$ we denote by $[a]_l$ the sum of its first $l$ terms.
Let us discuss the dimension of fibers of $\theta_n^{k,m,k-m}$.
Among the last $m$ coefficients of $j_1$ (which is a polynomial of
degree at most $k$) there are at most $\min (m, k-m+1)$ nonzero
ones. Knowing these coefficients and $\theta_n^{k,m,k-m}(j)$ one
can recover $j$. Thus fibers of $\theta_n^{k,m,k-m}$ are of
dimension $ \min (m, k-m+1)$.

For formal jets we define the similar morphism $ \hat{\theta}_n^m
: \hat{J}_{o_n}^m \C^n \to \hat{J}_E \tilde{\C}^n := \bigcup_{
\tilde{y} \in E} \hat{J}_{ \tilde{y} } \tilde{\C}^n$ given locally
by $( \hat{j}_1, \ldots , \hat{j}_n) \to (\hat{j}_1,
\hat{j}_2/\hat{j}_1, \ldots , \hat{j}_n/\hat{j}_1)$.

{\bf 4.5.} As usual we consider a closed embedding $(Y,y)
\hookrightarrow (\C^n , o_n)$. Let $ \tilde{Y}$ as a proper
transform of $Y$, $E_Y =E \cap \tilde{Y}$, and $\sigma_Y =
\sigma|_{\tilde{Y}}$, i.e. $\sigma_Y : \tilde{Y} \to Y$ is the
blowing-up of $(Y,y)$ at $y$. Put $\theta_Y^{k,m,l} =
\theta_n^{k,m,l}|_{ J_Y^{k,m}Y}$ and $ \hat{\theta}_Y^m
=\hat{\theta}_n^m |_{ \hat{J}_y^mY }$. It is easy to see that
$\theta_Y^{k,m,0} (J_y^{k,m}Y)\subset E_Y$.

\lemma {\it For every $ \tilde{y} \in E_Y$ the fiber $E^{
\tilde{y}}= (\theta_Y^{k,m,0}) ^{-1}( \tilde{y})$ is of dimension
at most $ \dim J_{ \tilde{y}}^{k-m} \tilde{Y} + \min (m, k-m+1)$
and for any $l =0, \ldots , k-m$ the image $ \theta_Y^{k,m,l}(E^{
\tilde{y} })$ is contained in $J_{ \tilde{y}}^l \tilde{Y}$.}

\proof The first statement follows from the second one for
$l=k-m$, the fact that $\theta_Y^{k,m,l} =\tau_Y^{r,l} \circ
\theta_Y^{k,m,r}$ for $r>l$, and the remark about the dimension of
$\theta_n^{k,m,k-m}$-fibers in 4.4. For the second statement put $
\tilde{j}= \theta_Y^{k,m,l}(j)$ where $j \in J_y^{k,m}Y$. Consider
a formal extension $ \hat{j} \in \hJ_yY$ of $j$, i.e. for every
regular function $h$ from the defining ideal of $(Y,y)$ in $(\C^n,
o_n)$ we have $h \circ \hat{j}=0$. Note that $\hat{\theta}_Y^m (
\hat{j})$ is a formal extension of $ \tilde{j}(t)$. Suppose that
$j_i$ and $\tilde{\xi}$ are as in 4.3 and 4.4, and
${\frac{d^m}{dt^m}} j_1(0) \ne 0$. Then for every regular function
$ \tilde{h}$ from the defining ideal of $( \tilde{Y}, \tilde{y})$
there exist $h$ as above and $s
>0$ so that $ h( \tilde{\xi}_1 , \tilde{\xi}_1 \tilde{\xi}_2,
\ldots , \tilde{\xi}_1 \tilde{\xi}_n) = \tilde{\xi}_1^s \tilde{h}
( \tilde{\xi}_1, \ldots , \tilde{\xi}_n)$. Hence $ \tilde{h} \circ
\hat{\theta}_Y^m ( \hat{j}) =0$ and $ \tilde{j } \in
J_{\tilde{y}}^{l} \tilde{Y}$. \qed

Induction on $k$ and Lemma 4.5 imply. \vs

\cor {\it The dimension of $J_{y}^{k,m}(Y)$ is at most $(k-m+1)
\dim Y + \min ( m-1, k-m)$. In particular, $\dim J_y^k (Y) \leq k
\dim Y$.} \vs

{\bf 4.6.} \lemma {\it Let $\varphi : (Y,y) \to (Z,z)$ be a
morphism, $j \in J_y^{k,m}(Y)$ be such that the multiplicity of $
\varphi^{(k)} (j)$ is $m$ (i.e. $\varphi^{(k)} (j) \in
J_z^{k,m}(Z)$). Let $\sigma: \tilde{Y} \to Y$ (resp. $ \delta
:\tilde{Z} \to Z$) be the blowing-up of $Y$ at $y$ (resp. $Z$ at
$z$) and $ \psi : \tilde{Y } - \to \tilde{Z}$ be the rational map
generated by $\varphi$. Then $\psi$ is regular at $ \tilde{y}=
\theta_Y^{k,m,0}(j)$ and sends it to $ \tilde{z}=
\theta_Z^{k,m,0}(\varphi^{(k)} (j))$. Furthermore, for $E^{
\tilde{y}}= (\theta_Y^{k,m,0}) ^{-1}( \tilde{y})$ we have
$\theta_Z^{k,m,l} \circ \varphi^{(k)}|_{E^{\tilde{y}}} =
\psi^{(l)} \circ \theta_Y^{k,m,l}|_{E^{\tilde{y}}}$.}

\proof Let $(Y,y) \hookrightarrow (\C^n , o_n)$ and $
\tilde{\C}^n$ be as in 4.5. In particular, $ \tilde{Y}$ can be
viewed as a subvariety of $ \tilde{\C}^n$. Let $(\C^s, o_s)$ and $
\tilde{\C}^s$ play the similar role for $ (Z,z)$. Then $\varphi$
is a restriction of a morphism $ \Phi : (\C^n, o_n) \to (\C^s,
o_s)$ which generates a rational map $ \Psi : \tilde{\C}^n - \to
\tilde{\C}^s$ such that $\psi$ is the restriction of $\Psi$. Thus
we can suppose that $(Y,y)=(\C^n, o_n)$ and $(Z,z)= (\C^s,o_s)$.
Let $ \tilde{\xi} =( \tilde{\xi}_1, \ldots, \tilde{\xi}_n)$ (resp.
$\tilde{\zeta} = ( \tilde{\zeta}_1, \ldots , \tilde{\zeta}_s)$ be
a local coordinate system on $ \tilde{\C}^n$ (resp. $
\tilde{\C}^s)$. Making linear coordinate changes we can suppose
that $ \tilde{y}= (1,0, \ldots ,0)$ in this local coordinate
system $ \tilde{\xi}$ (resp. $ \tilde{z} = (1,0, \ldots, 0)$ in $
\tilde{\zeta}$) and $\sigma ( \tilde{\xi}) = ( \tilde{\xi}_1,
\tilde{\xi}_2 \tilde{\xi}_1, \ldots \tilde{\xi}_n \tilde{\xi}_1 )$
(resp. $ \delta ( \tilde{\zeta} ) = ( \tilde{\zeta}_1,
\tilde{\zeta}_2 \tilde{\zeta}_1,  \ldots \tilde{\zeta}_s
\tilde{\zeta}_1)$). This implies that locally the coordinate form
of $\psi$ is $( \varphi_1 \circ \sigma, \varphi_2 \circ \sigma /
\varphi_1 \circ \sigma , \ldots \varphi_s \circ \sigma / \varphi_1
\circ \sigma)$ where $ \varphi =(\varphi_1, \ldots , \varphi_s)$,
and for every $j =(j_1, \ldots , j_n) \in E^{ \tilde{y}}$ we have
$ {\frac {d^m}{dt^m}} j_1(0) \ne 0$ (resp. $ {\frac {d^m} {dt^m}}
(\varphi_1 \circ j) (0) \ne 0$). Thus changing the coordinate $t$
on $(\C , 0)$ we can suppose that $j_1(t) =t^m$. Recall that for
every power series $a(t)$ the sum of its first $l$ terms is
denoted by $[a]_l$. Treating each $j_i$ as a polynomial of degree
at most $k$ we have $ \theta_s^{k,m,l} \circ \varphi^{(k)} (j)= (
[ \varphi_1 \circ j]_l, [\varphi_2 \circ j/ \varphi_1 \circ j]_l,
\ldots , [\varphi_s \circ j/ \varphi_1 \circ j]_l)$. Since $j_1(t)
=t^m$ we have on the other hand $ \sigma \circ \theta_n^{k,m,l}
(j) = ( [j_1]_l, [j_2/j_1]_l [j_1]_l, \ldots , [j_n/j_1]_l [j_1]_l
) = (j_1', [j_2]_{l+m}, \ldots , [j_n]_{l+m}) =: j'$ where $j_1'
=[j_1]_{l+m}$ for $l \geq m$, and $j_1'$ is zero for $l <m$. Hence
$\psi^{(l)} \circ \theta_n^{k,m,l} (j)= ( [ \varphi_1 \circ j']_l,
[\varphi_2 \circ j'/ \varphi_1 \circ j']_l, \ldots , [\varphi_s
\circ j'/ \varphi_1 \circ j']_l)$. As $ {\frac {d^m} {dt^m}}
(\varphi_1 \circ j) (0) \ne 0$ one can see that the last
expression coincides with those for $\theta_s^{k,m,l} \circ
\varphi^{(k)} (j)$. \qed

{\bf 4.7.} Let $h : Y_1 \to Y_2$ be a morphism of algebraic
varieties, $y_1 \in Y$, and $y_2=h(y_1)$. Then $h$ generates a
morphism $h_* : C_{y_1}Y_1 \to C_{y_2}Y_2$ of the reduced tangent
cones at $y_1$ and $y_2$ respectively where $h_*$ is just the
restriction of the induced linear map of the tangent spaces
$T_{y_1}Y_1 \to T_{y_2}Y_2$. It is known [D, Ch. 2.5.2] that if
$h$ is not unramified at $y_1$ (in particular, when it $y_1$ is
not a connected component of $h^{-1}(y_2)$) then the induced map
of (non-reduced) tangent cones is not an embedding. We need a
similar claim for reduced tangent cones.

\lemma {\it Let $h : (Y_1,y_1) \to (Y_2,y_2)$ be a morphism such
that $y_1 \ne h^{-1}(y_2) \cap Z_1$ for some irreducible analytic
branch $(Z_1,y_1)$ of $(Y_1,y_1)$. Let $(Z_2,y_2)$ be the proper
transform of $(Z_1,y_1)$ under $h$.

{\rm (1)} Then $h_*$ is not an embedding.

{\rm (2)} Let $V_i$ be the subspace of $T_{y_i}Y_i$ generated by
$C_{y_i}Z_i$ and $\dim V_1 \leq \dim V_2$. Then the closure of
$h_*(C_{y_1}Z_1)$ is a proper subvariety of $C_{y_2}Z_2$, i.e. $
\dim h_* (C_{y_1}Z_1) < \dim Y_2$.}

\proof Let $Y_i$ be a closed subvariety of $\C^n$ with coordinates
$x_1, \ldots , x_n$ so that $y_i$ is the origin. Consider the
homotety $(x_1, \ldots , x_n) \to (tx_1, \ldots , tx_n)$ where $t
\in \C^*$ and the image of $Y_i$ in $\C^n \simeq \C^n \times t$
under it. The closure of the union of these images is a subvariety
$ \check{Y}_i$ of $\C^{n+1}\simeq \C^n \times \C_t$ such that for
the natural projection $\tau_i : \check{Y}_i \to \C$ to the
$t$-axis, $\tau_i^{-1} (0)$ is isomorphic to ${C}_{y_i} Y_i$ and
there is an isomorphism $\varphi_i : Y_i \times \C^* \to
\tau_i^{-1}(\C^*)$ over $\C^*$ (e.g., see [D, Ch. 3.6.2]).
Moreover, $h$ generates a morphism $ \check{h} : \check{Y}_1 \to
\check{Y}_2$ such that $\tau_1 = \tau_2 \circ \check{h}$,
$\varphi_2^{-1} \circ \check{h} \circ \varphi_1 |_{Y_1 \times t}
=h$ for nonzero $t$, and $ \check{h}|_{\tau_1^{-1}(0)} =h_*$. The
closure $\check{Z}_i$ of $\varphi_i (Z_i \times \C^*)$ is an
irreducible analytic subvariety of $\check{Y}_i$ and $\check{Z}_i
\cap \tau_i^{-1}(0)$ is isomorphic to $C_{y_i}Z_i$. By the
assumption for any fixed $t \in \C^*$ the variety $ \hZ_1^t :=
\check{Z}_1 \cap \check{h}^{-1}(\varphi_2 (y_2 \times t)) \ne
\varphi_1 (y_1 \times t)$, i.e. $ \hZ_1^t$ is at least of
dimension 1 and the closure $\hZ_1$ of $\bigcup_{t \in \C^*}
\hZ_1^t$ is at least of dimension 2. Hence $\hZ_1^0 = \hZ_1 \cap
\tau_1^{-1}(0)$ contains a curve. Let $v_i$ be the vertex of
$C_{y_i}Y_i$. Note that $\varphi_i (y_i \times t)$ approaches
$v_i$ as $t \to 0$. Hence by continuity $ \check{h} (\hZ_1^0)
=v_2$, and, therefore, $h_*$ is not an embedding.

As $(Z_i,y_i)$ is an irreducible analytic branch of $(Y_i,y_i)$ it
follows easily from [M, Ch. 5A] that $C_{y_i}Z_i$ is irreducible.
Hence if the closure of $h_*(C_{y_1}Z_1)$ is not a proper
subvariety, it coincides with $C_{y_2}Z_2$ and, therefore,
$h_*(C_{y_1}Z_1)$ generates $V_2$. This implies that
$h_*(V_1)=V_2$ whence $h_*|_{V_1}: V_1 \to V_2$ is an isomorphism
as $\dim V_1 \leq \dim V_2$. Thus the restriction of $h_*$ to
$C_{y_1}Z_1 \subset V_1$ is an embedding contrary to (1). \qed

{\bf 4.8.} \lemma {\it Let $(Y_1,y_1) \to (Y_2,y_2) \to \ldots \to
(Y_s,y_s)$ be a sequence of birational morphisms of germs of
algebraic varieties and $g_{i_1,i_2}: (Y_{i_1},y_{i_2}) \to
(Y_{i_2},y_{i_2})$ be the composite morphisms for $i_1 < i_2$.
Suppose that $\dim Y_i=n\geq 2$ and $g_{i,i+1}^{-1}(y_{i+1}) \cap
Z_i \ne y_i$ for any irreducible analytic branch $(Z_i,y_i)$ of
$(Y_i,y_i)$. If $s \geq \dim T_{y_1}Y_1$ then $\dim
(g_{1,s})^{(1)}(C_{y_1}Y_1) \leq n-1$. }

\proof Let $V_i$ be the subspace of $T_{y_i}Y_i$  generated by
$C_{y_i}Z_i$ where $(Z_i,y_i)$ is the proper transform of
$(Z_1,y_1)$ under $g_{1,i}$, i.e. $g_{1,i}^{(1)}(V_1) \subset
V_i$. As $s \geq \dim T_{y_1}Y_1$ and $g_{1,i_2}=g_{i_1,i_2} \circ
g_{1,i_1}$ there exist $i_1 < i_2$ so that $\dim V_{i_1} \leq \dim
V_{i_2}$. By Lemma 4.7 $ \dim
g_{i_1,i_2}^{(1)}(C_{y_{i_1}}Y^1_{i_1}) < n$. As
$g_{1,i_1}^{(1)}(C_{y_1}Y_1^1) \subset C_{y_{i_1}}Y^1_{i_1}$ we
have the desired conclusion. \qed

{\bf 4.9.} In order to generalize the above Lemma to the case of
$k$-jet cones we need

\lemma {\it Let $(Y,y)$ be an irreducible germ of an analytic set,
$\sigma : \tilde{ Y} \to (Y,y)$ be its blowing-up at $y$, and
$E_Y$ be the exceptional divisor. Then $ E_Y \cap (\tilde{Z},
\tilde{y})$ is a divisor in $ \tilde{Y}$ where $( \tilde{Z},
\tilde{y})$ is any irreducible analytic branch of $\tilde{Y}$ at
any point $ \tilde{y} \in E_Y$.}

\proof Consider the union $U$ of all irreducible germs $(
\tilde{Z}, \tilde{y})$, $\tilde{y} \in E_Y$ of $ \tilde{Y}$ that
do not contain an open subset of $E_Y$. If $U \ne \emptyset$ then
$U$ is a proper analytic subset of $ \tilde{Y}$ of the same
dimension. Hence if $ {\nu} : \tY_{\nu} \to \tilde{Y}$ is
normalization then $ \tilde{Y}_{\nu}$ contains at least two
connected components: the proper transform $U'$ of $U$ and another
component $U''$ such that ${\nu} (U'') \supset E_Y$. There is a
natural proper morphism from $ \tilde{Y}_{\nu}$ into the
normalization $Y_{\nu}$ of $(Y,y)$. As $(Y,y)$ is irreducible the
preimage $y_{\nu}$ of $y$ in $Y_{\nu}$ is a point. But the
preimage of $y_{\nu}$ in $ \tilde{Y}_{\nu}$ is not connected (it
has points in both $U'$ and $U''$) in contradiction with the
Zariski Main Theorem. \qed

{\bf 4.10.} \lemma {\it Let the assumption of Lemma 4.8 hold and
$k>0$. Suppose that $s = (2l) ^{k-1}l$ and $l \geq \max_Z \dim TZ
$ where $Z$ is the result of any sequence of $r$ blowing-ups of
any $(Y_i,y_i)$ at $y_i$ and infinitely near points with $1\leq i
\leq s$ and $0 \leq r \leq k-1$. Then $\dim g_{1,s}^{(k)}
(J_{y_1}^{k,m} Y_1) \leq (k-m+1) (n-1)+ \min (m-1, k-m)$ for every
$1 \leq m \leq k$. In particular, $ \dim g_{1,s}^{(k)} (
J_{y_1}^kY_1 ) \leq k (n-1)$.}

\proof Let $k=m$. Then every $j \in J_{y_i}^{k,k}Y_i$ is of form
$j=j_0 \circ h$ where $j_0 \in C_{y_i}Y_i$ and $h : \C \to \C , t
\to t^k$. Hence $g_{i_1,i_2}^{(k)}(j) = g_{i_1,i_2}^{(1)}(j_0)
\circ h$. In this case the statement follows from Lemma 4.8. In
particular, Lemma is true for $k=1$. We use now induction on $k$
and inside it induction on $k-m$. Let $s_0=s/2+1$ and for $i<s_0$
let $S_i^0$ be the subvariety of $J_{y_i}^{k,m}Y_i$ such that
$g_{i,s_0}^{(k)}(S_i^0) \subset J_{y_{s_0}}^{k,m+1}Y_{s_0}$. By
induction $g_{1,s}^{(k)}(S_1^0) \subset
g_{s_0,s}^{(k)}(J_{y_{s_0}}^{k,m+1}Y_{s_0})$ is of dimension at
most $(k-m)(n-1) + \min (m, k-m-1)$. Thus (since $n>1$) it
suffices to consider jets from $S_1$ where $S_i =J_{y_i}^{k,m}Y_i
\setminus S_i^0$.

Let  $ \tilde{Y}_i$ be the blowing-up of $Y_i$ at $y_i$. Its
exceptional divisor $E_i$ is naturally isomorphic to the base of
the cone $C_{y_i}Y_i$ and $g_{i_1,i_2}$ generates a birational map
$h_{i_1,i_2}: \tilde{Y}_{i_1} - \to \tilde{Y}_{i_2}$.  Deleting
the indeterminacy points (i.e. replacing $E_i$ by its Zariski open
subset $E_i^*$) we can suppose that $h_{i_1,i_2}$ is regular on $
\tilde{Y}_{i_1}^*= (\tilde{Y}_i \setminus E_i) \cup E_i^*$ for
$i_2 \leq s_0$. Note that for every $j \in S_{i_1}$ the
multiplicity of $g_{i_1,i_2}(j)$ is $m$ whence by Lemma 4.6
$\theta_{Y_i}^{k,m,0}(S_i) \subset E_i^*$. By Lemma 4.8 for every
$(2l)^{k-1}
> q\geq 0$ there exist $l q < i_1 < i_2 \leq l (q+1)$ such that
the dimension of $g_{i_1,i_2}^{(1)}(C_{y_{i_1}}Y_{i_1})$ is at
most $n-1$ whence $\dim h_{i_1,i_2}(E_{i_1}^*) \leq n-2$, i.e.
$E_{i_1}^*$ is the exceptional divisor of $h_{i_1,i_2}$. Put
$e_{i_1,i_2}=h_{li_1,li_2}$ and $Z_i =\tilde{Y}_{li}^*$. We get a
sequence of birational morphisms $Z_1 \to Z_2 \to \ldots \to
Z_{s_1}$ where $s_1=(2l) ^{k-2}l$. Note that $E_{li_1}$ is an
exceptional divisor of $e_{i_1,i_2}$ and, by Lemma 4.9 it meets
every irreducible analytic branch $(Z_i^1,z_i)$ of $(Z_i,z_i)$
where $z_i$ is any point of $E_{li_1}^*$. Thus this new sequence
of birational morphisms satisfies the assumption of this Lemma. By
induction $\dim e_{1,s_1}^{(k-m)}(J_{z_1}^{k-m}Z_1) \leq (k-m)
(n-1)$, and Lemma 4.6 implies that $\theta_{Y_{ls_1}}^{k,m,k-m}
\circ g_{1, ls_1 }^{(k)}(S_1) \subset \bigcup_{ \tilde{y}_1 \in
E^*_1} h_{1,ls_1}^{(k-m)} (J_{ \tilde{y}_1}^{k-m}\tilde{Y}_1)
\subset \bigcup_{z_1 \in h_{1,l}(E^*_1)} e_{1,s_1}^{(k-m)}
(J_{z_1}^{k-m}Z_1)$. As $\dim h_{1,l}(E^*_1) \leq n-2$, we have $
\dim \theta_{Y_{ls_1}}^{k,m,k-m} \circ g_{1, ls_1 }^{(k)}(S_1)
\leq (k-m+1)(n-1)-1$. Taking into consideration the remark about
the dimension of $\theta_n^{k,m,k-m}$-fibers in 4.4 we get the
desired conclusion. \qed

{\bf 5. The proof of Theorems B and A.\\ 5.1} By Lemma 2.1 we can
suppose that $X$ and $Y$ are normal in Theorem B. In the case when
$n= \dim Y=1$ the result follows from the fact that a bijective
morphism of smooth curves is an isomorphism. Consider $n>1$.
Suppose that under the assumption of Theorem B (i) $f$ is an
automorphism and $g$ is not. By Corollary 2.4 there exists an
exceptional divisor $D$ for $g$. Assume $x$ be a general point in
$\rho^{-1}(D) \ne \emptyset$. In particular, $x_s =f^s (x)$ is
also a general point in $\rho^{-1}(D)$. Let $\psi : Y \to \C^n$ be
a dominant morphism. As $f^s$ is an automorphism $\dim (\psi \circ
\rho \circ f^s) ^{(k)} (J_x^kX) = ( \psi \circ \rho ) ^{(k)}
(J_{x_s}^kX) = ( \psi \circ \rho ) ^{(k)} (J_{x}^kX)$ as both $x$
and $x_s$ are general points. By Proposition 3.2 there exists
$n_0$ such that for any $k$ we have $\dim ( \psi \circ \rho )
^{(k)} (J_{x}^kX) \geq k n - n_0$. On the other hand $Y$ is
locally analytically irreducible since it is normal. By Lemma 4.10
$ \dim (\psi \circ g^s \circ \rho) ^{(k)}(J_x^kX) \leq k (n-1)$
for sufficiently large $s$. Since $\psi \circ \rho \circ f^s =
\psi \circ g^s \circ \rho$ we get a contradiction which proves
that $g$ is an automorphism. The proof of Theorem B (ii) is
similar. \qed

{\bf 5.2.} Theorem B yields Theorem A in the case of $k = \C$. We
need to reduce the general case to this one. Let $\bar {k}$ be the
algebraic closure of the field $k$. Recall that $\bar {k}$ is a
faithfully flat $k$-module.
This means that an endomorphism $\varphi : T \to T$ of a
$k$-algebra $T$ is an automorphism iff the endomorphism $\varphi
\otimes_k Id_{\bar {k}} : T \otimes_k \bar {k} \to T \otimes_k
\bar {k}$ is an automorphism. Thus we can replace the rings $S$
and $R$ in Theorem A by $S \otimes_k \bar {k}$ and $R \otimes_k
\bar {k}$ respectively. That is, we can suppose from the beginning
that $k$ is algebraically closed. We consider case (i) only since
the other case is similar. It is equivalent to the analogue of
Theorem B (i) in which $X$ and $Y$ are already affine algebraic
varieties over $k$. Note that Lemmas 2.1 and 2.2 hold for every
algebraically closed field whence we can suppose that $Y$ is
normal and $g$ is birational. Hence if we assume that $g$ is not
an automorphism then a coordinate function of $g^{-1}$ has a pole
at a point $y_0 \in Y$. Let $k'$ be the subfield of $k$ generated
by a finite number of elements which include the coordinates of
$y_0$ in the ambient Euclidean space, the coefficients of
coordinate functions of $\rho , g, f$ and $f^{-1}$ (as polynomials
over $k$), and the coefficients of generators of the defining
ideals of $X$ and $Y$. Consider our varieties and morphisms over
$k'$ instead of $k$ and denote the corresponding objects by
$X',Y',f',g'$, and $\rho'$. Note that $g'$ is not an automorphism
as $y_0 \in Y'$.
But $k'$ can be embedded as a subfield in $\C$ by the ``Lefschetz
principle" [BCW]. Hence theorem B (i) implies that the coordinate
functions of $(g')^{-1}$ cannot have a pole at $y_0$.
Contradiction. \qed

\begin{center} REFERENCES \end{center} \vs
\noin [A] J. Ax, {\em Injective endomorphisms of varieties and
schemes}, Pacific J. Math., {\bf 31} (1969), 1--7. \\ \noin [BCW]
H. Bass, E. Connell, D. Wright, {\em The Jacobian conjecture:
reduction of degree and the formal expansion of the
inverse}, Bull. Amer. Math. Soc. {\bf 7} (1982), 287--330. \\
\noin [BN] S. Bell, R. Narasimhan, {\em Proper holomorphic
mappings of complex spaces}, Encyclopaedia of Math. Sci., Several
Compl. Var., {\bf 69}, Springer-Verlag, Berlin-Heidelberg-New
York, e.a., 1990, 1--38.\\ \noin [CZ] E. Connell, J. Zweibel, {\em
Subrings invariant under polynomial maps}, Houston Journal of
Math., {\bf 20} (1994), 175--185. \\ \noin [D] V. I. Danilov, {\em
Algebraic varieties and schemes}, Algebraic Geometry I
(Encyclopaedia of Math. Sci., v. {\bf 23}) (Russian), VINITI,
Moscow, 1988.  \\
\noin [H] R. Hartshorne, {\em Algebraic geometry},
Springer-Verlag, New York-Berlin-Heidelberg, 1977.\\ \noin [I] S.
Iitaka {\em On logarithmic Kodaira dimension for algebraic
varieties}, In: Complex Analysis and Algebraic Geometry, 175--189,
Kinokuniya, Tokyo, 1978. \\
\noin [M-L] L. Makar-Limanov, {\em On the hypersurface $x + x^2y +
z^2 + t^3 = 0$ in ${\C}^{4}$ or a ${\C}^3$-like threefold which is
not ${\C}^3$}, Israel J. Math. {\bf 96} (1996), 419--429.\\ \noin
\noin [M] D. Mumford, {\em Algebraic geometry}, Springer-Verlag,
Berlin-Heidelberg- New York, e.a., 1976. \\ \noin [P] Th.
Peternell, {\em Modifications}, Encyclopaedia of Math. Sci.,
Several Compl. Var., {\bf 74}, Springer-Verlag, Berlin-Heidelberg-
New York, e.a.,
1994, 221--258. \\
\end{document}